\documentclass[a4paper,dvipsnames]{amsart}

\usepackage[ngerman,english]{babel}
\usepackage[TS1,OT1,T1]{fontenc}
\usepackage[latin1]{inputenc}
\usepackage{amsmath}
\usepackage{amsthm}
\usepackage{mathbbol}
\usepackage{amssymb}         %% diverse Matheerweiterungen, z.B. \mathbb{R}
\usepackage{epsfig}          %% um eps-Dateien einzubinden (\epsfig{file=...})
\usepackage{longtable}       %% fuer Tabellen ueber mehrere Seiten
\usepackage{multicol}        %% fuer mehrspaltigen Text auf einer Seite

\usepackage{graphicx}               %% einbinden von Grafiken
\usepackage{setspace}
\usepackage{colortbl}
\usepackage{lscape}
\usepackage{fancybox}
\usepackage[arrow, matrix, curve]{xy}

\usepackage{mathtools}
\usepackage[perpage]{footmisc}

\usepackage[a4paper,top=3cm,bottom=2cm,left=3cm,right=3cm,marginparwidth=1.75cm]{geometry}

\renewcommand{\thefootnote}{\fnsymbol{footnote}}

\renewcommand*\thefootnote{%
  \ifcase\value{footnote}\or
  *\or
  **\or
  ***\or
  ****\or
  *****\or
  ******\or
  *******\or
  ********\or
  *********\fi
}

\mathtoolsset{showonlyrefs} \mathtoolsset{showmanualtags}

\addtocontents{toc}{\protect\thispagestyle{empty}}
\usepackage{babelbib}

\makeindex

\theoremstyle{plain}
\newtheorem{Sa}{Theorem}
\newtheorem{Le}{Lemma}
\newtheorem*{Kor}{Corollary}

\theoremstyle{definition}
\newtheorem{Def}{Definition}

\theoremstyle{remark}
\newtheorem*{Beh}{Proposition}
\newtheorem{Bsp}{Example}
\newtheorem*{Bem}{Note}

\newenvironment{acknowledgements}{\textbf{Acknowledgements}\\}{}

\newcommand{\id}{\operatorname{id}}
\newcommand{\dd}{\mathrm{d}}

\newcommand{\Bild}{\operatorname{im}}
\newcommand{\Kern}{\operatorname{kern}}
\newcommand{\Lin}{\operatorname{Lin}}

\newcommand{\gklammer}[1]{\left\{#1\right\}}
\newcommand{\klammer}[1]{\left(#1\right)}
\newcommand{\eklammer}[1]{\left[#1\right]}

\newcommand{\ds}[1]{$\displaystyle{#1}$}
\newcommand{\re}{\operatorname{Re}}
\newcommand{\im}{\operatorname{Im}}

\newcommand{\BIGOP}[1]{\mathop{\mathchoice%
{\raise-0.22em\hbox{\huge $#1$}}%
{\raise-0.05em\hbox{\Large $#1$}}{\hbox{\large $#1$}}{#1}}}

\newcommand{\BIGboxplus}{\mathop{\mathchoice%
{\raise-0.35em\hbox{\huge $\boxplus$}}%
{\raise-0.15em\hbox{\Large $\boxplus$}}{\hbox{\large
$\boxplus$}}{\boxplus}}}
\usepackage[dvipsnames]{xcolor}
\usepackage[pagebackref,hypertexnames=false, colorlinks, citecolor=BrickRed,linkcolor=MidnightBlue, urlcolor=BrickRed]{hyperref}
\usepackage[hyphenbreaks]{breakurl}

\title{Additive Deformations of Hopf Algebras}
\author{Malte Gerhold}
\address{Universität Greifswald\\ Institut für Mathematik und Informatik\\
Walther-Rathenau-Straße 47\\
17487 Greifswald, Germany}
\thanks{\ \vspace*{-1em}}
\thanks{Final version published at \url{https://doi.org/10.1016/j.jalgebra.2011.05.019}}
\thanks{\textcopyright 2023. This manuscript version is made available under the CC-BY-NC-ND 4.0 license\\ \url{https://creativecommons.org/licenses/by-nc-nd/4.0/}.}

\begin{document}

\begin{abstract}
Additive deformations of bialgebras in the sense of Wirth are deformations of the multiplication map of the bialgebra fulfilling a compatibility condition with the coalgebra structure and a continuity condition. Two problems concerning additive deformations are considered.

With a deformation theory a cohomology theory should be developed. Here a variant of the Hochschild cohomology is used. The main result in the first partad of this paper is the characterization of the trivial deformations, i.e.\ deformations generated by a coboundary.

Starting with a Hopf algebra, one would expect the deformed multiplications to have some analogue of the antipode, which we call deformed antipodes. We prove, that deformed antipodes always exist, explore their properties, give a formula to calculate them given the deformation and the antipode of the original Hopf algebra and show in the cocommutative case, that each deformation splits into a trivial part and into a part with constant antipodes.
\end{abstract}

\maketitle

% \begin{keyword}
% Hopf algebras\sep deformation theory\sep quantum algebras\sep quantum groups\sep cohomology
% \end{keyword}

\section{Introduction}\label{sec:intro}

Deformations of algebras are closely related to cohomology as Gerstenhaber showed in his papers \cite{Gers63} and \cite{Gers64}. Suppose that $\mathcal{A}$ is an algebra and $(\mu_t)_{t\geq0}$ a family of associative multiplications on $\mathcal{A}$, which can in any sense be written in the form
\begin{equation}
	\mu_t(a\otimes b)=\mu(a\otimes b)+tF(a\otimes b)+\mathcal{O}(t^2),
\end{equation}
where $\mu_0=\mu$ is the original multiplication of the algebra. Writing down the associativity condition for $\mu_t$ and comparing the terms of first order yields that
\begin{equation}
	\mu (F(a\otimes b)\otimes c)+F(\mu(a\otimes b)\otimes c)=\mu (a\otimes F(b\otimes c))+F(a\otimes\mu (b\otimes c))
\end{equation}
and after rearranging
\begin{equation}
	aF(b\otimes c)-F(ab\otimes c)+F(a\otimes bc)-F(a\otimes b)c=0,
\end{equation}
so the infinitesimal deformation $F$ is a cocycle in the Hochschild cohomology associated with the $\mathcal{A}$-Bimodule structure on $\mathcal{A}$ given by multiplication.

Additive deformations were first introduced by Wirth in \cite{Wir02} as a generalization of Weyl algebras. Given a finite dimensional complex vector space $\mathcal{V}$ with an alternating bilinear form $\mathfrak{s}$ (if $\mathfrak{s}$ is nondegenerate, this is a symplectic form, whence the letter $\mathfrak{s}$) one can form the Weyl algebra $H_\mathfrak{s}:=T(\mathcal{V})/I_\mathfrak{s}$, where $T(\mathcal{V})=\bigoplus_{n=0}^\infty \mathcal{V}^{\otimes n}$ is the tensor algebra over $\mathcal{V}$ and $I_\mathfrak{s}$ is the ideal generated by elements of the form $v\otimes w-w\otimes v-\mathfrak{s}(v,w)$, so that $\mathfrak{s}$ becomes the commutator in the algebra $H_\mathfrak{s}$. It is clear, that $H_0=:H$ is the symmetric tensor algebra over $\mathcal{V}$ and it can be shown that the family $(H_{t\mathfrak{s}})_{t\in\mathbb{R}}$ can be identified with a deformation of the symmetric algebra, i.e.\ there are invertible linear mappings $\Phi_t:H_{t\mathfrak{s}}\rightarrow H$ and we get a family $(\mu_t)_{t\in\mathbb{R}}$ of multiplcations on $H$ (see \cite{Wir02} and references therein). Setting in particular $\mathcal{V}$ the vectorspace with basis $\gklammer{a,a^\dagger}$ and $\mathfrak{s}(a,a^\dagger)=\hbar$ the obtained algebra is the algebra of the quantum harmonic oscillator. In \cite{Maj95} Majid just calls this a bialgebra like structure. It is in fact an example of an additive deformation in the following sense.

An additive deformation of a bialgebra $\mathcal{B}$ is a familiy $(\mu_t)_{t\in\mathbb{R}}$ of multiplications, which turns $\mathcal{B}_t=(\mathcal{B},\mu_t,\Eins)$ into a unital algebra ($\Eins$ is the unit element of the original algebra $\mathcal{B}$) such that $\Delta:\mathcal{B}_{t+s}\rightarrow\mathcal{B}_t\otimes\mathcal{B}_s$ is an algebra homomorphism and which satisfies some continuity condition (see Definition\,\ref{Def:additive_deformation}). Wirth proved in \cite{Wir02} that all additive deformations are of the form $\mu_t=\mu\star e_\star^{tL}$, where $L$ is a commuting normalized 2-cocycle in the Hochschild cohomology associated with the $\mathcal{B}$-bimodule structure on $\mathbb{C}$ given by the counit (see Section\,\ref{sec:cohomology} and Theorem\,\ref{Sa:generator_of_an_additive_deformation}). 

Wirth also showed in \cite{Wir02}, that a Schoenberg correspondence holds for additive deformations. In \cite{Sch05} and \cite{G09} quantum Lévy processes on additive deformations are introduced, so additive deformations are of interest in quantum probability.

In the present paper we have two goals. First we wish to introduce a cohomology, such that we have a 1-1-correspondence between additive deformations and all cocycles. This also gives a concept of trivial deformations, i.e.\ deformations generated by a coboundary. We give a handy characterization of these trivial deformations. The second goal is to describe additive deformations of Hopf algebras. When one starts with a Hopf algebra, one would expect the deformed multiplications to have some analogue to the antipode, which we call deformed antipodes. We prove the existence of such deformed antipodes and describe their behaviour.

In Section\,\ref{sec:cohomology} we introduce a cohomology, such that the generators of additive deformations are exactly the 2-cocycles. This is done by requiring each $n$-cochain to be normalized and to commute with $\mu^{(n)}$, the multiplication map for $n$ factors. One has to show that this is a cochain complex, explicitly, that $\partial c$ is normalized and commuting if $c$ is. The same can be done for $*$-deformations of $*$-algebras.

Once the cohomology is established the question is, what kind of deformations are generated by coboundaries. It is shown, that those deformations are of the form
\begin{equation}
\mu_t=\Phi_t\circ\mu\circ (\Phi_t^{-1}\otimes\Phi_t^{-1})
\end{equation}
where the $\Phi_t$ constitute a pointwise continuous one parameter group of invertible linear operators on $\mathcal{B}$ that commute in the sense that 
\begin{equation}
(\Phi_t\otimes\id )\circ\Delta= (\id\otimes\Phi_t )\circ\Delta.
\end{equation}
When $L=\partial\psi$ is the generator of the additive deformation, then $\Phi_t=(\id\otimes e_\star^{-t\psi})\circ\Delta$ is the one parameter group of operators.

The second section of the paper discusses additive deformations of Hopf algebras. Deforming the multiplication of a bialgebra $\mathcal{B}$ also gives a deformed convolution product $\star_t$ for linear maps from $\mathcal{B}$ to $\mathcal{B}$
\begin{equation}
A\star_t B:=\mu_t\circ(A\otimes B)\circ\Delta,
\end{equation}
where $(\mu_t)_{t\in\mathbb{R}}$ is a deformation of the multiplication map $\mu$ of $\mathcal{B}$.
If $\mathcal{B}$ is a Hopf algebra, i.e.\ there is an antipode $S$, which is the convolution inverse of the identity on $\mathcal{B}$ w.r.t.\ $\star=\mu\circ(\cdot\otimes\cdot)\circ\Delta$, it is quite natural to ask, whether there are also deformed antipodes $S_t$, which fulfill
\begin{equation}\label{intro:eq:deformed_antipodes1}
\mu_t\circ (S_t\otimes\id )\circ\Delta= \mu_t\circ (\id\otimes S_t)\circ\Delta 
=\Eins\delta
\end{equation}
i.e.\ they are convolution inverse to the identity map w.r.t.\ $\star_t$ and if they exist, what properties they have. Such a deformation is called a Hopf deformation.

In a Hopf algebra the antipode $S$ is an algebra antihomomorphism and a coalgebra antihomomorphism, i.e.
\begin{align}
S\circ\mu &=\mu\circ (S\otimes S)\circ\tau,\\
\Delta\circ S&=\tau\circ(S\otimes S)\circ\Delta.
\end{align}
Similar properties hold for the deformed antipodes $S_t$ of a Hopf deformation. We can prove
\begin{align}
S_t\circ\mu_{-t} &=\mu_t\circ\tau \circ (S_t\otimes S_t),\label{intro:eq:algebra_anti_hom}\\
\Delta\circ S_{t+r}&=(S_t\otimes S_r)\circ\tau\circ\Delta.\label{intro:eq:koalgebra_anti_hom}
\end{align}

Applying $\delta\otimes\delta$ to \eqref{intro:eq:koalgebra_anti_hom} we get
\begin{equation}
\delta\circ S_{t+r}=((\delta\circ S_t)\otimes (\delta\circ S_r))\circ\tau\circ\Delta
=((\delta\circ S_r)\otimes (\delta\circ S_t))\circ\Delta,
\end{equation}
i.e.\ $\delta\circ S_t$ is a convolution semigroup w.r.t.\ $\star=(\cdot\otimes\cdot)\circ\Delta$. So one would like to prove that this semigroup has a generator, such that the $S_t$ are of the form
\begin{equation}\label{intro:eq:deformed_antipodes2}
S_t=S\star e_\star^{-t\sigma}.
\end{equation}

To get a hint, how to find $\sigma$, we assume for the moment that $\delta\circ S_t$ is differentiable in $0$ and define
\begin{equation}\label{intro:eq:sigma1}
\sigma:=-\left.\frac{\dd}{\dd t}\delta\circ S_t\right|_{t=0}.
\end{equation}
Then we can apply $\delta$ to \eqref{intro:eq:deformed_antipodes1} and differentiate to get
\begin{equation}
L\circ(S\otimes\id)\circ\Delta-\sigma=L\circ(\id\otimes S)\circ\Delta-\sigma=0
\end{equation}
or after rearranging
\begin{equation}\label{intro:eq:sigma2}
\sigma=L\circ(S\otimes\id)\circ\Delta=L\circ(\id\otimes S)\circ\Delta.
\end{equation}
In fact we will prove that every additive deformation of a Hopf algebra is a Hopf deformation and \eqref{intro:eq:deformed_antipodes2} and \eqref{intro:eq:sigma2} give a formula for the deformed antipodes.

In two special cases the structure can even be better understood. In the case of a trivial deformation it is easy to see that
\begin{equation}
S_t=\Phi_t\circ S\circ\Phi_t
\end{equation}
is another way to find the deformed antipodes. Differentiating this also gives a second formula for the generator
\begin{equation}
\sigma=\psi+\psi\circ S.
\end{equation}
If the bialgebra $\mathcal{B}$ is cocommutative, we show that every additive deformation splits in a trivial part and a part with constant antipodes. 
Applying $\delta$ to \eqref{intro:eq:algebra_anti_hom} and differentiating yields
\begin{equation}
-\sigma\circ\mu-L=L\circ(S\otimes S)\circ\tau-\sigma\otimes\delta-\delta\otimes\sigma
\end{equation}
or after rearranging
\begin{equation}\label{intro:eq:del_sigma}
L+L\circ (S\otimes S)\circ\tau=\delta\otimes\sigma-\sigma\circ\mu+\sigma\otimes\delta=\partial\sigma.
\end{equation}
So $L$ can be written as
\begin{equation}
L=\underbrace{\frac{1}{2}\partial\sigma}_{:=L_1}
+\underbrace{\frac{1}{2}(L-L\circ(S\otimes S)\circ\tau)}_{:=L_2}
\end{equation}
and if $\mathcal{B}$ is cocommutative the second part corresponds to constant antipodes.

\section{Notation and Basic Definitions}\label{sec:basics}

All vector spaces considered are over the complex numbers, denoted by $\mathbb{C}$. The algebraic dual of a vector space $\mathcal{V}$ is denoted $\mathcal{V}':=\gklammer{\varphi:\mathcal{V}\rightarrow\mathbb{C}~|~\text{$\varphi$ linear}}$.  The tensor product $\otimes$ is the usual algebraic tensor product. If $\mathcal{V}$ is a vector space we write 
\begin{equation}
\mathcal{V}^{\otimes n}:=\underbrace{\mathcal{V}\otimes\dots\otimes\mathcal{V}}_{n\times}
\end{equation}
for $n\geq1$ and $\mathcal{V}^{\otimes 0}:=\mathbb{C}$.

A bialgebra $(\mathcal{B},\mu,\Eins,\Delta,\delta)$ is a complex unital associative algebra $(\mathcal{B},\mu,\Eins)$ for which the mappings $\Delta:\mathcal{B}\rightarrow\mathcal{B}\otimes\mathcal{B}$ and $\delta:\mathcal{B}\rightarrow\mathbb{C}$ are algebra homomorphisms and satisfy coassociativity and counit property respectively. A Hopf algebra is a bialgebra with an antipode, i.e.\ a linear mapping $S:\mathcal{B}\rightarrow\mathcal{B}$ with
\begin{equation}
\mu\circ (\id\otimes S)\circ\Delta=\Eins\delta=\mu\circ (S\otimes\id )\circ\Delta.
\end{equation}
A $*$-bialgebra is a bialgebra with an involution, i.e.\ $(\mathcal{B},\mu,\Eins,*)$ is a $*$-algebra and $\Delta,\delta$ are $*$-homomorphisms. If $\mathcal{B}$ is a $*$-bialgebra, an involution on $\mathcal{B}\otimes\mathcal{B}$ is given by $(a\otimes b)^*=a^*\otimes b^*$. A $*$-Hopf algebra is a Hopf algebra which also is a $*$-bialgebra. For details on Hopf algebras and bialgebras see e.g.\ \cite{Abe80, Swe69} for $*$-Hopf algebras \cite{KlSc97}.

We use Sweedler's notation, writing $\Delta a=\sum_{k=0}^n a^{(1)}_k\otimes a^{(2)}_k=:a_{(1)}\otimes a_{(2)}$ and the notations $\mu^{(n)}:\mathcal{B}^{\otimes n}\rightarrow\mathcal{B}$, $\Delta^{(n)}:\mathcal{B}\rightarrow\mathcal{B}^{\otimes n}$
\begin{align}
\mu^{(0)}(\lambda)&=\lambda\Eins&
\Delta^{(0)}&=\delta\\
\mu^{(n+1)}&=\mu\circ(\id\otimes\mu^{(n)})&
\Delta^{(n+1)}&=(\id\otimes\Delta^{(n)})\circ\Delta.
\end{align}
The Sweedler notation for this is
\begin{equation}
\Delta^{(n)}a=a_{(1)}\otimes\dots\otimes a_{(n)}.
\end{equation}

% The counit applied to an $a_{(k)}$ can be just left out due to the counit property and if $S(a_{(k)})$ is multiplied with $a_{(k+1)}$ or $a_{(k-1)}$ both are just left out and the remaining indices are adapted, for example
% \begin{equation}
% a_{(1)}\delta(a_{(2)})S(a_{(3)})\otimes a_{(4)}=a_{(1)}S(a_{(2)})\otimes a_{(3)}=\Eins\otimes a.
% \end{equation}

With $\mathcal{B}$ also each $\mathcal{B}^{\otimes n}$ is a bialgebra in the natural way. We frequently use the comultiplication on $\mathcal{B}\otimes\mathcal{B}$, which we denote by $\Lambda$ and which is defined by
\begin{equation}
\Lambda(a\otimes b)=a_{(1)}\otimes b_{(1)}\otimes a_{(2)}\otimes b_{(2)},
\end{equation}
i.e.\ $\Lambda=(\id\otimes\tau\otimes\id)\circ(\Delta\otimes\Delta)$. The counit of $\mathcal{B}\otimes\mathcal{B}$ is just $\delta\otimes\delta$.

If $(\mathcal{C},\Delta)$ is a coalgebra and $(\mathcal{A},m)$ is an algebra, we define the convolution product for mappings $R,S:\mathcal{C}\rightarrow\mathcal{A}$ by $R\star S:=m\circ(R\otimes S)\circ\Delta$. In our context $\mathcal{C}$ and $\mathcal{A}$ are usually tensor powers of the same bialgebra $\mathcal{B}$.

A pointwise continuous convolution semigroup is a family $(\varphi_t)_{t\geq0}$ of linear mappings $\varphi_t:\mathcal{B}\rightarrow\mathbb{C}$ such that
\begin{itemize}
\item $\varphi_t\star\varphi_s=\varphi_{t+s}$
\item $\varphi_t(b)\xrightarrow{t\rightarrow0}\delta(b)\quad\forall b\in\mathcal{B}$
\end{itemize}
Note that $\delta$ is the unit for the multiplication $\star$ on $\mathcal{B}'$ (This is exactly the counit property).
It follows from the fundamental theorem for coalgebras, that for a pointwise continuous convolution semigroup there exists a generator $\psi$, which is the pointwise limit
\begin{equation}
\psi(b)=-\left.\frac{\dd \varphi_t(b)}{\dd t}\right|_{t=0}=-\lim_{t\rightarrow0}\frac{\varphi_t(b)-\delta(b)}{t}
\end{equation}
and for which we have
\begin{equation}
\varphi_t=e_\star^{-t\psi}.
\end{equation}
Cf.\ \cite{ASvW88} section 4 for details.

\begin{Def}\label{Def:additive_deformation}
An \emph{additive deformation} of the bialgebra $\mathcal{B}$ is a family $(\mu_t)_{t\geq0}$ of mappings $\mu_t:\mathcal{B}\otimes\mathcal{B}\rightarrow\mathcal{B}$ such that
\begin{enumerate}
	\item $(\mathcal{B},\mu_t,\Eins)$ is a unital algebra for each $t\geq0$
	\item $\mu_0=\mu$
	\item $\Delta\circ\mu_{t+s}=(\mu_t\otimes\mu_s)\circ(\id\otimes\tau\otimes\id)\circ(\Delta\otimes\Delta)$ ($\tau$ denotes the flip on $\mathcal{B}\otimes\mathcal{B}$)
	\item the mapping $t\mapsto \delta\circ\mu_t$ is pointwise continuous, i.e.\ $\delta\circ\mu_t\xrightarrow{t\rightarrow0}\delta\circ\mu
=\delta\otimes\delta$ pointwise
	\item if $\mathcal{B}$ is a $*$-bialgebra and for each $t\geq0$ $(\mathcal{B},\mu_t,\Eins,*)$ is a unital $*$-algebra, we call the deformation an additive deformation of a $*$-bialgebra.
\end{enumerate}
\end{Def}
\goodbreak

The following theorem was first proven by Wirth in \cite{Wir02}. A proof can also be found in \cite{G09}.

\begin{Sa}\label{Sa:generator_of_an_additive_deformation}
Let $(\mu_t)_{t\geq0}$ be an additive Deformation of the bialgebra $\mathcal{B}$. Then $L=\left.\frac{\dd(\delta\circ\mu_t)}{\dd t}\right|$ exists pointwise and we have for $a,b,c\in \mathcal{B},t\geq0$
\begin{enumerate}
\item\ds{\mu_t=\mu\star e_{\star}^{tL}}
\item\ds{\mu\star L=L\star\mu\quad\text{'$L$ is commuting'}}
\item\ds{L(\Eins\otimes\Eins)=0\quad\text{'$L$ is normalized'}}
\item\ds{\delta(a)L(b\otimes c)-L(ab\otimes c)+L(a\otimes bc)-L(a\otimes b)\delta(c)=0}.\\
'$L$ is a coboundary'
\item[] \hspace{-20pt} If $(\mu_t)_{t\geq0}$ is a $*$-bialgebra deformation, then
\item\ds{L(b\otimes c)=\overline{L(c^*\otimes b^*)}}\quad\text{'$L$ is hermitian'}
\end{enumerate}
also holds.

Conversely, if $L:\mathcal{B}\otimes\mathcal{B}\rightarrow\mathbb{C}$ is a linear mapping, which fulfills conditions 2,3 and 4 (in case of $*$-bialgebra also 5), than the first equation defines an additive deformation on $\mathcal{B}$.
\end{Sa}

\section{Cohomology of Additive Deformations}\label{sec:cohomology}

\subsection{Subcohomologies of the Hochschild cohomology}\label{subsec:hom}

A cochain complex consists of a sequence of vector spaces $C=(C_n)_{n\in\mathbb{N}}$ and linear mappings $\partial_n:C_n\rightarrow C_{n+1}$ such that $\partial_{n+1}\circ\partial_n=0$ for all $n\in\mathbb{N}$. The elements of $Z_n(C)=\Kern\partial_n$ are called ($n-$)cocycles, the elements of $B_n(C)=\Bild\partial_{n-1}$ are called ($n-$)coboundaries and the vector-space $H_n(C)=Z_n(C)/B_n(C)$ is called $n$-th cohomology. A sequence $D=(D_n)_{n\in\mathbb{N}}$ is called subcomplex, if $D_n\subseteq C_n$ and $\partial_n D_n\subseteq D_{n+1}$ for all $n$. Then $\klammer{(D_n)_{n\in\mathbb{N}},(\left.\partial_n\right|_{D_n})_{n\in\mathbb{N}}}$ is again a cochain complex and we have:
\begin{enumerate}
	\item The cocycles of $D$ are exactly the cocycles of $C$, belonging to $D$, i.e.
	\begin{equation}
		Z_n(D)=Z_n(C)\cap D_n,
	\end{equation}
	\item each coboundary of $D$ is a coboundary of $C$, i.e.
	\begin{equation}
		B_n(D)\subseteq B_n(C)\cap D_n,
	\end{equation}
	\item equality holds, iff the mapping $H_n(D)\rightarrow H_n(C),f+B_n(D)\mapsto f+B_n(C)$ is an injection,
	\item If $D,E$ are subcomplexes, then $(D_n\cap E_n)_{n\in\mathbb{N}}$ is a subcomplex.
\end{enumerate}
Points 1,2 and 4 are obvious, while 3 follows from the observation, that the kernel of the given mapping is exactly $B_n(C)\cap D_n$.

For an algebra $\mathcal{A}$ and an $\mathcal{A}$-bimodule $M$ we define
\begin{equation}
	C_n:=\Lin(\mathcal{A}^{\otimes n},M)=\gklammer{f:\mathcal{A}^{\otimes n}\rightarrow M~|~f\text{ linear}}.
\end{equation}
One can show, that together with the coboundary operator
\begin{multline}
	\partial_n f (a_1,\dots,a_{n+a}):=\\
		a_1.f(a_2,\dots,a_{n+1})+\sum_{i=1}^n (-1)^i~f(a_1,\dots,a_ia_{i+1},\dots,a_{n+1})
        +(-1)^{n+1}~f(a_1,\dots,a_n).a_{n+1}
\end{multline}
this is a cochain complex. Especially for $\mathcal{A}=\mathcal{B}$ a bialgebra and $M=\mathbb{C}$ the $\mathcal{B}$-bimodule given by $a.\lambda.b=\delta(a)\lambda\delta(b)$ for $\lambda\in\mathbb{C}$ and $a,b\in\mathcal{B}$ we have
\begin{multline}\label{eq:cobop}
	\partial_n f (a_1,\dots,a_{n+a}):=\delta(a_1)f(a_2,\dots,a_{n+1})+\sum_{i=1}^n (-1)^i~f(a_1,\dots,a_ia_{i+1},\dots,a_{n+1})\\
        +(-1)^{n+1}~f(a_1,\dots,a_n)\delta(a_{n+1}).
\end{multline}

The generators of additive deformations are normalized commuting cocycles, so it is natural to define
\begin{align}
	C_n^{(\mathbf{N})}&=\gklammer{f\in C_n~|~f(\Eins^{\otimes n})=0},\\
	C_n^{(\mathbf{C})}&=\gklammer{f\in C_n~|~f\star\mu^{(n)}=\mu^{(n)}\star f}.
\end{align}
If $\mathcal{B}$ is a $*$-bialgebra the generators are also hermitian. We define for $f\in C_n$
\begin{equation}
\widetilde{f}(a_1\otimes\dots\otimes a_n):=\overline{f(a_n^*\otimes\dots\otimes a_1^*)}
\end{equation}
and set
\begin{equation}
C_n^{(\mathbf{H})}=
\begin{cases}
\gklammer{f\in C_n~\Big|~\widetilde{f}=f},\textcolor{white}{\bigg|}&\text{ if $\left\lceil\frac{n}{2}\right\rceil$ odd, i.e.\ $n=1,2,5,6,\dots$}\\
\gklammer{f\in C_n~\Big|~\widetilde{f}=-f},&\text{ if $\left\lceil\frac{n}{2}\right\rceil$ even, i.e.\ $n=0,3,4,7,8,\dots$}
\end{cases}
\end{equation}

\begin{Beh}
$C_n^{(\mathbf{N})}$, $C_n^{(\mathbf{C})}$ and $C_n^{(\mathbf{H})}$ are subcomplexes of $C_n$.
\end{Beh}

\begin{proof}
We only need to show that $\partial C_n^{(*)}\subseteq C_n^{(*)}$ for $*=\mathbf{N,C,H}$.
\begin{description}
\item[$\mathbf{N}$:]
Let $f\in C_n^{(\mathbf{N})}$, i.e.\ $f(\Eins^{\otimes n})=0$. Then
\begin{equation}
\partial f(\Eins^{\otimes (n+1)})=
\delta(\Eins)f(\Eins^{\otimes n})+\sum_{i=1}^n (-1)^i~f(\Eins^{\otimes n})
        +(-1)^{n+1}~f(\Eins^{\otimes n})\delta(\Eins)=0
\end{equation}
\item[$\mathbf{C}$:]
For $f\in C_n^{(\mathbf{C})}$ we get
\begin{multline}
\partial f\star\mu^{(n+1)}=\\
\klammer{\delta\otimes f+\sum_{k=1}^n (-1)^k
f\circ(\id_{k-1}\otimes\mu\otimes\id_{n-k})
+(-1)^{n+1}f\otimes\delta}
\star\mu^{(n+1)}.
\end{multline}
Next we show, that each summand commutes with $\mu$ under convolution:
\begin{align}
&\quad(\delta\otimes f)\star\mu^{(n+1)}(a_1\otimes\dots\otimes
        a_{n+1})\\
            &=\delta(a_1^{(1)})f\klammer{a_2^{(1)}\otimes\dots\otimes
                a_{n+1}^{(1)}}a_1^{(2)}\dots a_{n+1}^{(2)}\\
            &=a_1 f\klammer{a_2^{(1)}\otimes\dots\otimes
                a_{n+1}^{(1)}}a_2^{(2)}\dots a_{n+1}^{(2)}\\
            &=a_1 f\klammer{a_2^{(2)}\otimes\dots\otimes
                a_{n+1}^{(2)}}a_2^{(1)}\dots a_{n+1}^{(1)}
                \qquad\quad\text{(as $f\star\mu^{(n)}=\mu^{(n)}\star f$)}\\
            &=a_1^{(1)}\dots a_{n+1}^{(1)}\delta(a_1^{(2)})
                f\klammer{a_2^{(2)}\otimes\dots\otimes a_{n+1}^{(2)}}\\
            &=\mu^{(n+1)}\star(\delta\otimes f)(a_1\otimes\dots\otimes
        a_{n+1}).
\end{align}
Analoguesly we see that
\begin{equation}
(f\otimes\delta)\star\mu^{(n+1)}=\mu^{(n+1)}\star(f\otimes\delta).
\end{equation}
For the remaining summands we calculate
\begin{align}
        &\quad(f\circ(\id_{k-1}\otimes\mu\otimes\id_{n-k}))\star\mu^{(n+1)}
            (a_1\otimes\dots\otimes a_{n+1})\\
        &=f\klammer{a_1^{(1)}\otimes\dots\otimes
            (a_k^{(1)}a_{k+1}^{(1)})\otimes\dots\otimes
            a_{n+1}^{(1)}}
            a_1^{(2)}\dots a_k^{(2)}a_{k+1}^{(2)}\dots
            a_{n+1}^{(2)}\\
        &=f\klammer{a_1^{(1)}\otimes\dots\otimes
            (a_k a_{k+1})^{(1)}\otimes\dots\otimes a_{n+1}^{(1)}}
            a_1^{(2)}\dots (a_ka_{k+1})^{(2)}\dots
            a_{n+1}^{(2)}
        \intertext{
            \raggedleft{(as $\Delta$ is an algebra-homomorphism)}
            }\displaybreak[0]\\
        &=f\klammer{a_1^{(2)}\otimes\dots\otimes
            (a_k a_{k+1})^{(2)}\otimes\dots\otimes a_{n+1}^{(2)}}
            a_1^{(1)}\dots (a_ka_{k+1})^{(1)}\dots
            a_{n+1}^{(1)}
        \intertext{
            \raggedleft{(as $f\star\mu^{(n)}=\mu^{(n)}\star f$)}
            }
        &=\mu^{(n+1)}\star(f\circ(\id_{k-1}\otimes\mu\otimes\id_{n-k}))
            (a_1\otimes\dots\otimes a_{n+1}).
\end{align}
\item[$\mathbf{H}$:]
Let $\widetilde{f}=\pm f$. For $n$ odd, we get
\begin{align}
&\quad\widetilde{\partial f}(a_1,\dots,a_{n+1})=\overline{\partial f(a_{n+1}^*,\dots,a_1^*)}
=\overline{\delta(a_{n+1}^*)f(a_n^*,\dots,a_{1}^*)}\\
&\quad
+\sum_{i=1}^n (-1)^{n+1-i}~\overline{f(a_{n+1}^*,\dots,a_{i+1}^*a_i^*,\dots,a_{1}^*)}
+\overline{f(a_{n+1}^*,\dots,a_2^*)\delta(a_{1}^*)}\\
&=\delta(a_1)\widetilde{f}(a_2,\dots,a_{n+1})+\sum_{i=1}^n (-1)^i~\widetilde{f}(a_1,\dots,a_ia_{i+1},\dots,a_{n+1})+\widetilde{f}(a_1,\dots,a_n)\delta(a_{n+1})\\
&=\pm\partial f(a_1,\dots,a_{n+1})
\end{align}
and for $n$ even
\begin{align}
&\quad\widetilde{\partial f}(a_1,\dots,a_{n+1})
=\overline{\partial f(a_{n+1}^*,\dots,a_1^*)}
=\overline{\delta(a_{n+1}^*)f(a_n^*,\dots,a_{1}^*)}\\
&\quad+\sum_{i=1}^n (-1)^{n+1-i}~\overline{f(a_{n+1}^*,\dots,a_{i+1}^*a_i^*,\dots,a_{1}^*)}-\overline{f(a_{n+1}^*,\dots,a_2^*)\delta(a_{1}^*)}\\
&=-\delta(a_1)\widetilde{f}(a_2,\dots,a_{n+1})-\sum_{i=1}^n (-1)^i~\widetilde{f}(a_1,\dots,a_ia_{i+1},\dots,a_{n+1})+\widetilde{f}(a_1,\dots,a_n)\delta(a_{n+1})\\
&=\mp\partial f(a_1,\dots,a_{n+1}).
\end{align}
\end{description}
\end{proof}
Since the intersection of subcomplexes is again a subcomplex we have
\begin{Kor}
$C_n^{(\mathbf{NC})}:=C_n^{(\mathbf{N})}\cap C_n^{(\mathbf{C})}$ and 
$C_n^{(\mathbf{NCH})}:=C_n^{(\mathbf{NC})}\cap C_n^{(\mathbf{H})}$ are cochain complexes with the coboundary operator \eqref{eq:cobop}.
\end{Kor}

\subsection{Characterization of the trivial deformations}\label{subsec:def}

For an additive deformation of the bialgebra $\mathcal{B}$ the generator $L$ of the convolution-semigroup $(\delta\circ\mu_t)_{t\geq0}$ is an element of $Z_2^{(\mathbf{NC})}$ and conversely if $L\in Z_2^{(\mathbf{NC})}$ we can define an additive deformation via $\mu_t:=\mu\star e_\star^{tL}$. In the case of a $*$-bialgebra the generators are exactly the elements of $Z_2^{(\mathbf{NCH})}$. We wish to answer the question which deformations are generated by the coboundaries, i.e.\ the elements of $B_2^{(\mathbf{NC})}$ or $B_2^{(\mathbf{NCH})}$ respectively.

\begin{Sa}\label{Sa:triv_def}
Let $\mathcal{B}$ be a bialgebra, $L\in B_2^{(\mathbf{NC})}$, $L=\partial\psi$ with $\psi\in C_1^{(\mathbf{NC})}$. Then $(\Phi_t)_{t\geq0}$ is a pointwise continuous semigroup of unital algebra isomorphisms
$\Phi_t:(\mathcal{B},\mu)\rightarrow(\mathcal{B},\mu_t)$, for which
\begin{equation}\label{eq:comm}
(\Phi_t\otimes\id)\circ\Delta=(\id\otimes\Phi_t)\circ\Delta\quad\text{for all $t\geq0$},
\end{equation}
where $\Phi_t=\id\star e_\star^{-t\psi}$ and $\mu_t=\mu\star e_\star^{tL}$.
If $\mathcal{B}$ is a $*$-bialgebra and $L\in B_2^{(\mathbf{NCH})}$, then we can choose $\psi\in C_1^{(\mathbf{NCH})}$ and the $\Phi_t$ are $*$-algebra isomorphisms.

Conversely, if $(\Phi_t)_{t\geq0}$ is a pointwise continuous semigroup of invertible linear mappings $\Phi_t:\mathcal{B}\rightarrow\mathcal{B}$, such that $\Phi_t(\Eins)=\Eins$ for all $t\geq0$, and which fulfills \eqref{eq:comm}, then 
\begin{equation}
\mu_t:=\Phi_t\circ\mu\circ(\Phi_t^{-1}\otimes\Phi_t^{-1})
\end{equation}
defines an additive Deformation of $\mathcal{B}$ with generator $L\in B_2^{(\mathbf{NC})}$.
If $\mathcal{B}$ is a $*$-algebra and the $\Phi_t$ are hermitian, then we get an additive deformation of a $*$-bialgebra and $L\in B_2^{(\mathbf{NCH})}$.
\end{Sa}

Before we prove this, we recall the following Lemma.
When $\mathcal{B}$ is a bialgebra and $\varphi:\mathcal{B}\rightarrow\mathbb{C}$ a linear functional on $\mathcal{B}$ we define
\begin{equation}
R_\varphi:\mathcal{B}\rightarrow\mathcal{B},\quad
R_\varphi:=\id\star\varphi=(\id\otimes\varphi)\circ\Delta.
\end{equation}

\begin{Le}
For $\varphi,\psi\in\mathcal{B}'$ the following holds:
\begin{enumerate}
\item \ds{R_{\varphi}\circ R_{\psi}=R_{\varphi\star\psi}}
\item \ds{\delta\circ R_{\varphi}=\varphi}
\item \ds{R_{\delta\otimes\varphi}=\id\otimes R_\varphi}
\item \ds{R_{\varphi\otimes\delta}=R_\varphi\otimes\id}
\item \ds{\mu\circ R_{\varphi\circ\mu}=R_\varphi\circ\mu}
\end{enumerate}
% \begin{align}
% R_{\varphi}\circ R_{\psi}&=R_{\varphi\star\psi}\\
% \delta\circ R_{\varphi}&=\varphi\\
% R_{\delta\otimes\varphi}&=\id\otimes R_\varphi\\
% R_{\varphi\otimes\delta}&=R_\varphi\otimes\id\\
% \mu\circ R_{\varphi\circ\mu}&=R_\varphi\circ\mu
% \end{align}
Note that the last three equations are between operators on the bialgebra $\mathcal{B}\otimes\mathcal{B}$.
\end{Le}

\begin{proof}
This is all straightforward to verify.
\end{proof}

\begin{proof}[Proof of Theorem\,\ref{Sa:triv_def}.]
Let $\mathcal{B}$ be a bialgebra and $L\in B_2^{(\mathbf{NC})}$ a coboundary, $L=\partial\psi$ with $\psi\in C_1^{(\mathbf{NC})}$. We write $\varphi_t:=e_\star^{-t\psi}$ and note, that this is a pointwise continuous convolution semigroup and the $\varphi_t$ are commuting (i.e.\ $\varphi_t\star\id=\id\star\varphi_t$) since $\psi$ is. Then the mappings $\Phi_t=R_{\varphi_t}$ yield a pointwise continuous semigroup of linear operators on $\mathcal{B}$ and we only need to show, that they are unital algebra isomorphisms. It is obvious, that $\Phi_t(\Eins)=\Eins$, since $\psi(\Eins)=0$, and $\Phi_t\circ\Phi_{-t}=\id$, so $\Phi_t$ is invertible. We have to prove that $\Phi_t:(\mathcal{B},\mu)\rightarrow(\mathcal{B},\mu_t)$ is an algebra homomorphism, i.e.\  
\begin{equation}
\mu_t=\mu\star e_\star^{tL}=\Phi_t\circ\mu\circ(\Phi_t^{-1}\otimes\Phi_t^{-1}).
\end{equation}
From
\begin{align}
e_\star^{tL}&=e_\star^{t\partial\psi}\\
&=e_\star^{t(\delta\otimes\psi-\psi\circ\mu+\psi\otimes\delta)}\\
&=e_\star^{-t\psi\circ\mu}\star e_\star^{t\delta\otimes\psi}\star e_\star^{t\psi\otimes\delta}\\
&=(\varphi_{t}\circ\mu)\star(\delta\otimes\varphi_{-t})\star(\varphi_{-t}\otimes\delta),
\end{align}
where we used that $\delta\otimes\psi$, $\psi\circ\mu$ and $\psi\otimes\delta$ commute under convolution, we conclude
\begin{align}
\mu_t&=\mu\star e_\star^{t\partial\psi}\\
&=(\mu\otimes e_\star^{t\partial\psi})\circ\Lambda\\
&=\mu\circ R_{e_\star^{t\partial\psi}}\\
&=\mu\circ R_{(\varphi_{t}\circ\mu)\star(\delta\otimes\varphi_{-t})\star(\varphi_{-t}\otimes\delta)}\\
&=\mu\circ R_{\varphi_{t}\circ\mu}\circ (\id\otimes R_{\varphi_{-t}})\circ(R_{\varphi_{-t}}\otimes\id)\\
&=R_{\varphi_{t}}\circ\mu\circ(R_{\varphi_{-t}}\otimes R_{\varphi_{-t}})\\
&=\Phi_t\circ\mu\circ(\Phi_t^{-1}\otimes\Phi_t^{-1}).
\end{align}
It is clear that the $\Phi_t$ are $*$-homomorphisms in the $*$-bialgebra case.

Now let $(\Phi_t)_{t\geq0}$ be pointwise continuous semigroup of invertible linear mappings with $\Phi_t(\Eins)=\Eins$ and $(\Phi_t\otimes\id)\circ\Delta=(\id\otimes\Phi_t)\circ\Delta$. Then we write $\varphi_t=\delta\circ\Phi_t$ and observe that
\begin{enumerate}
\item $(\varphi_t)_{t\geq0}$ is a pointwise continuous convolution semigroup, so that there is a $\psi\in C_1^{(\mathbf{NC})}$ with $\varphi_t=e_\star^{-t\psi}$. Indeed
\begin{align}
\varphi_t\star\varphi_s&=((\delta\circ\Phi_t)\otimes(\delta\circ\Phi_s))\circ\Delta\\
&=(\delta\otimes\delta)\circ(\Phi_t\otimes\id)\circ(\id\otimes\Phi_s)\circ\Delta\\
&=(\delta\otimes\delta)\circ(\id\otimes\Phi_t)\circ(\id\otimes\Phi_s)\circ\Delta\\
&=(\delta\otimes\delta)\circ(\id\otimes\Phi_{t+s})\circ\Delta\\
&=\varphi_{t+s}
\end{align}
and $\psi(\Eins)=0,\psi\star\id=\id\star\psi$ follow from $\varphi_t(\Eins)=1$ and $\varphi_t\star\id=\id\star\varphi_t$ via differentiation. If the $\Phi_t$ are hermitian, $\psi$ is also hermitian, i.e.\ $\psi\in C_1^{(\mathbf{NCH})}$.
\item $\Phi_t=R_{\varphi_t}$, as
\begin{align}
R_{\varphi_t}&=(\id\otimes(\delta\circ\Phi_t))\circ\Delta\\
&=(\id\otimes\delta)\circ(\id\otimes\Phi_t)\circ\Delta\\
&=(\id\otimes\delta)\circ(\Phi_t\otimes\id)\circ\Delta\\
&=\Phi_t.
\end{align}
\end{enumerate}
So the first part of the theorem tells us, that $L=\partial\psi\in B_2^{(\mathbf{NC})}$ is the generator of an additive deformation, for which
\begin{equation}
\mu_t=\Phi_t\circ\mu\circ(\Phi_t^{-1}\otimes\Phi_t^{-1}).
\end{equation}
If $\mathcal{B}$ is a $*$-bialgebra and all the $\Phi_t$ are hermitian, then so are all the $\varphi_t$ and via differentiation also $\psi$. That means $L\in B_2^{(\mathbf{NCH})}$ and the deformation is an additive deformation of a $*$-bialgebra.
\end{proof}

\begin{Bem}
  Let $L\in B_2$, i.e.\ $L=\partial\psi$ for a linear functional
  $\psi$. Then it follows that $\psi(\Eins)=0$ iff
  $L(\Eins\otimes\Eins)=0$ and if $L$ is hermitian $1/2
  (\psi+\widetilde{\psi})$ is a hermitian functional whose coboundary
  is $L$. In other words $B_2^{(\mathbf{N})}=B_2\cap
  C_2^{(\mathbf{N})}$ and $B_2^{(\mathbf{H})}=B_2\cap
  C_2^{(\mathbf{H})}$, but it is not clear under which circumstances
  $B_2^{(\mathbf{C})}=B_2\cap C_2^{(\mathbf{C})}$ holds, i.e.\ if
  there are 2-coboundaries that commute but are not coboundaries of a
  commuting functional. This possible difference is actually the main
  reason why we need the altered cochain complex to get a good notion
  of trivial deformations.
\end{Bem}

\section{Additive Deformations of Hopf Algebras}\label{sec:hopf}

\subsection{Definition of Hopf deformations and general observations}

\begin{Le}
Let $\mathcal{B}$ be a Bialgebra and $L$ generator of an additive deformation. Then we can define
\begin{equation}
\mu_t:=e_\star^{tL}\star\mu
\end{equation}
for all $t\in\mathbb{R}$ (i.e.\ not only for $t\geq0$) and
\begin{equation}
\Delta:\mathcal{B}_{t+s}\rightarrow\mathcal{B}_t\otimes\mathcal{B}_s
\end{equation}
is an algebra homomorphism for all $s,t\in\mathbb{R}$.
\end{Le}

\begin{proof}
It follows from Theorem\,\ref{Sa:generator_of_an_additive_deformation} that $-L$ is the generator of an additive deformation, so for $t<0$ the definition of $\mu_t$ yields a multiplication on $\mathcal{B}$. We calculate
\begin{equation}\begin{split}
\Delta\circ\mu_{t+s}
&=\Delta\circ(\mu\otimes e_\star^{(t+s)L})\circ\Lambda\\
&=((\Delta\circ\mu)\otimes e_\star^{(t+s)L})\circ\Lambda\\
&=(\mu\otimes\mu\otimes e_\star^{tL}\otimes e_\star^{sL})\circ\Lambda^{(4)}\\
&=(\mu\otimes e_\star^{tL}\otimes\mu\otimes e_\star^{sL})\circ\Lambda^{(4)}\\
&=((\mu\star e_\star^{tL})\otimes(\mu\star e_\star^{sL}))\circ\Lambda\\
&=(\mu_t\otimes\mu_s)\circ\Lambda
\end{split}\end{equation}
\end{proof}

From now on we always view an additive deformation as a family of multiplications indexed by all real numbers.

\begin{Def}\label{Def:Hopf_deformation}
An additive deformation is called a \emph{Hopf deformation}, if for all $t\in\mathbb{R}$ there exists a linear mapping $S_t:\mathcal{B}\rightarrow\mathcal{B}$ such that
\begin{equation}\label{eq:deformed_antipode}
\mu_t\circ(S_t\otimes\id)\circ\Delta=\mu_t\circ(\id\otimes S_t)\circ\Delta=\Eins\delta.
\end{equation}
\end{Def}

For $t=0$ this of course implies, that $\mathcal{B}$ is a Hopf algebra with antipode $S=S_0$.

\begin{Bem}
Many proofs in this section follow a common path. To show an identity $a=b$, we find an element $c$ and a convolution product $\diamond$ such that $a\diamond c=c\diamond b=\delta$ where $\delta$ is the neutral element for $\diamond$. Then we conclude
\begin{equation}
a=a\diamond\delta=a\diamond c\diamond b=\delta\diamond b=b
\end{equation}
and hence
\begin{equation}
a=b=c^{-1}.
\end{equation}
\end{Bem}

Let $\mathcal{B}$ be a bialgebra with additive deformation $(\mu_t)_{t\in\mathbb{R}}$ and $\mu_t=\mu\star e_\star^{tL}$ for a normalized, commuting cocycle $L$.

\begin{Sa}
If a family $S_t$ with \eqref{eq:deformed_antipode} exists, it is uniquely determined and the following statements hold:
\begin{enumerate}
\item $S_t(\Eins)=\Eins$
\item $S_t:\mathcal{B}_{-t}\rightarrow \mathcal{B}_t$ is an algebra antihomomorphism, i.e.
\begin{equation}\label{eq:S_t_mu}
S_t\circ\mu_{-t}=\mu_t\circ (S_t\otimes S_t)\circ\tau
\end{equation}\label{eq:S_t_Delta}
\item $\displaystyle{\Delta\circ S_{t+r}=(S_t\otimes S_r)\circ\tau\circ\Delta}$

\item If $\mathcal{B}$ is cocommutative, i.e.\ $\Delta=\tau\circ\Delta$, then $S_t$ is invertible
for all $t\in\mathbb{R}$ and $(S_t)^{-1}=S_{-t}$.
\end{enumerate}
\end{Sa}

\begin{proof}
\begin{description}
\item[(Uniqueness)]{
The uniqueness statement is clear, as \eqref{eq:deformed_antipode} states, that $S_t$ is the two-sided convolution inverse of the identity mapping on $\mathcal{B}$ w.r.t.\ $\star_t$.}
\item[1.]{This is clear, since 
\begin{equation}
\Eins=\mu_t\circ(S_t\otimes\id)\circ\Delta(\Eins)=S_t(\Eins).
\end{equation}
}
\item[2.]{We show, that both sides of \eqref{eq:S_t_mu} are convolution inverses of $\mu_t$ w.r.t.\ $\star_t$:
\begin{equation}\begin{split}
(S_t\circ\mu_{-t})\star_t\mu_t&=\mu_t\circ(S_t\otimes\id)\circ(\mu_{-t}\otimes\mu_t)\circ\Lambda\\
&=\mu_t\circ(S_t\otimes\id)\circ\Delta\circ\mu=\delta\circ\mu\Eins=\delta\otimes\Eins\delta
\end{split}\end{equation}
and
\begin{equation}\begin{split}
&\quad\mu_t\star_t(\mu_t\circ (S_t\otimes S_t)\circ\tau)\, (a\otimes b)\\
&=\mu_t\circ(\mu_t\otimes\mu_t)\circ (\id_2\otimes ((S_t\otimes S_t)\circ\tau ))
\circ\Lambda\, (a\otimes b)\\
&=\mu_t^{(4)}(a_{(1)}\otimes b_{(1)}\otimes S_t(b_{(2)})\otimes S_t(a_{(2)}))\\
&=\delta(b)\mu_t (a_{(1)}\otimes S_t(a_{(1)}))\\
&=\delta(a)\delta(b)\Eins.
\end{split}\end{equation}}
\item[3.]{For linear maps from the coalgebra $(B,\Delta)$ to the algebra $(B_t\otimes B_r)$ we have a convolution $\diamond$ defined as
\begin{equation}
A\diamond B=(\mu_t\otimes\mu_r )\circ(\id\otimes\tau\otimes\id)\circ (A\otimes B)\circ\Delta.
\end{equation}
We show that both sides of \eqref{eq:S_t_Delta} are inverses of $\Delta$ w.r.t.\ $\diamond$:
\begin{equation}\begin{split}
(\Delta\circ S_{t+r})\diamond\Delta
&=(\mu_t\otimes\mu_r )\circ(\id\otimes\tau\otimes\id)\circ (\Delta\otimes\Delta)\circ
(S_{t+r}\otimes\id)\circ\Delta\\
&=\Delta\circ\mu_{t+r}\circ(S_{t+r}\otimes\id)\circ\Delta=\delta\,\Delta(\Eins)=\Eins\delta\otimes\Eins
\end{split}\end{equation}
and
\begin{equation}\begin{split}
&\quad\Delta\diamond( (S_t\otimes S_r)\circ\tau\circ\Delta) (a)\\
&=(\mu_t\otimes\mu_r )\circ(\id\otimes\tau\otimes\id)\circ
	(\id_2\otimes S_t\otimes S_r)\circ(\id_2\otimes\tau)\circ\Delta^{(4)}(a)\\
&=(\mu_t\otimes\mu_r )(a_{(1)}\otimes S_t(a_{(4)})\otimes a_{(2)}\otimes S_r(a_{(3)}))\\
&=\mu_t (a_{(1)}\otimes S_t(a_{(2)}))\otimes\Eins\\
&=\delta(a)\Eins\otimes\Eins
\end{split}\end{equation}
}
\item[4.]{
Let $\Delta=\tau\circ\Delta$. Then
\begin{equation}\begin{split}
(S_t\circ S_{-t})\star_t S_t
&=\mu_t\circ(S_t\otimes S_t)\circ(S_{-t}\otimes\id)\circ\Delta\\
&=S_t\circ\mu_{-t}\circ\tau\circ(S_{-t}\otimes\id)\circ\Delta\\
&=S_t\circ\mu_{-t}\circ(\id\otimes S_{-t})\circ\tau\circ\Delta\\
&=\delta S_t(\Eins)=\Eins\delta
\end{split}\end{equation}
}
\end{description}
\end{proof}

\subsection{The deformed antipodes for trivial deformations}

\begin{Sa}
Let $\mathcal{B}$ be a Hopf algebra and $(\mu_t)_{t\in\mathbb{R}}$ a trivial deformation,
\begin{equation}
\mu_t=\Phi_t\circ\mu\circ(\Phi_t^{-1}\otimes\Phi_t^{-1}),
\end{equation}
and
\begin{equation}
\Phi_t=\id\star e_\star^{-t\psi}
\end{equation}
for a commuting, normalized linear functional $\psi$.
Then 
\begin{equation}
S_t=\Phi_t\circ S\circ\Phi_t= S\star e_\star^{-t(\psi\circ S+\psi)} 
\end{equation}
is the deformed antipode, so $(\mu_t)_{t\in\mathbb{R}}$ is a Hopf deformation.
\end{Sa}

\begin{proof}
All we have to show is that 
\begin{equation}
\mu_t\circ(S_t\otimes\id)\circ\Delta=\mu_t\circ(\id\otimes S_t)\circ\Delta=\Eins\delta,
\end{equation}
for $S_t=\Phi_t\circ S\circ\Phi_t$ and $S_t=S\star e_\star^{-t(\psi\circ S+\psi)}$.
In the case $S_t=\Phi_t\circ S\circ\Phi_t$ we calculate
\begin{equation}\begin{split}
\mu_t\circ(S_t\otimes\id)\circ\Delta
&=\Phi_t\circ\mu\circ (\Phi_t^{-1}\otimes\Phi_t^{-1})
	\circ ((\Phi_t\circ S\circ\Phi_t )\otimes\id )\circ\Delta\\
&=\Phi_t\circ\mu\circ (S\circ\Phi_t\otimes \Phi_t^{-1})\circ\Delta\\
&=\Phi_t\circ\mu\circ (S\otimes \Phi_t^{-1})\circ(\Phi_t\otimes\id)\circ\Delta\\
&=\Phi_t\circ\mu\circ (S\otimes \Phi_t^{-1})\circ(\id\otimes\Phi_t)\circ\Delta\\
&=\Phi_t\circ\mu\circ (S\otimes \id)\circ\Delta\\
&=\delta\Phi_t(\Eins)=\Eins\delta
\end{split}\end{equation}
and the second equality is proven in the same way.

Now we consider the case $S_t=S\star e_\star^{-t(\psi\circ S+\psi)}$. We first recall that $\psi$ ist commuting and $L=\partial\psi$ is the generator of the additive deformation. Next we observe that
\begin{equation}\begin{split}
(\psi\circ S)\star S &= (\psi\otimes\id)\circ (S\otimes S)\circ\Delta\\
&=(\psi\otimes\id )\circ\tau\circ\Delta\circ S\\
&=(\psi\otimes\id )\circ\Delta\circ S\\
&=(\id\otimes\psi )\circ (S\otimes S)\circ\Delta\\
&=S\star(\psi\circ S).
\end{split}\end{equation}
With this in mind we calculate
\begin{equation}\begin{split}
\mu_t\circ(S_t\otimes\id)\circ\Delta\,(a)&=
(\mu\otimes e_\star^{tL})\circ\Lambda 
(e_\star^{-t\psi} (S(a_{(1)})) e_\star^{-t\psi}(a_{(2)})S(a_{(3)})\otimes a_{(4)})\\
&=e_\star^{-t\psi}(S(a_{(1)}))e_\star^{-t\psi}(a_{(2)})e_\star^{tL}(S(a_{(3)})\otimes a_{(4)})\\
&=\delta(a)
\end{split}\end{equation}
since
\begin{equation}\begin{split}
&\quad e_\star^{tL}(S(a_{(1)})\otimes a_{(2)})\\
&=e_{\star}^{t\delta\otimes\psi}(S(a_{(1)})\otimes a_{(2)})
e_\star^{-t\psi\circ\mu}(S(a_{(3)})\otimes a_{(4)})
e_{\star}^{t\psi\otimes\delta}(S(a_{(5)})\otimes a_{(6)})\\
&=e_{\star}^{t\psi}(a_{(1)})e_{\star}^{t\psi}(S(a_{(2)})).
\end{split}\end{equation}
Again the second equality is proven similarly.

One can also prove this by writing $\Phi_t=(e_\star^{-t\psi}\otimes\id)\circ\Delta$ in $S_t=\Phi_t\circ S\circ\Phi_t$ and using that $S,\psi$ and $\psi\circ S$ all commute with each other.
\end{proof}

It is still possible that the deformed antipodes are constant. We have

\begin{Sa}
Let $L$ be generator of a trivial additive deformation. Then it has constant antipodes, i.e.\ $S_t=S\ \forall t\in\mathbb{R}$ iff
\begin{equation}
\Phi_t\circ S=S\circ\Phi_{-t}.
\end{equation}
for all $t\in\mathbb{R}$.
\end{Sa}

\begin{proof}
This follow directly from $S_t=\Phi_t\circ S\circ\Phi_t$ and $\Phi_t^{-1}=\Phi_{-t}$.
\end{proof}

\subsection{The deformed antipodes in the general case}

We want to show, that every additive deformation of a Hopf algebra is
a Hopf deformation and give a formula for the deformed antipodes.

\begin{Le}\label{Le:sigma}
We have
\begin{equation}
L\circ(S\otimes\id)\circ\Delta = L\circ(\id\otimes S)\circ\Delta.
\end{equation}
\end{Le}

\begin{proof}
  From $\partial L=0$ it follows easily that $L(a\otimes\Eins)=L(\Eins\otimes
  a)=0$ for all $a\in\mathcal{B}$. Hence
\begin{equation}
\begin{split}
0	&=\partial L(a_{(1)}\otimes S(a_{(2)})\otimes a_{(3)}) \\
	&=\delta(a_{(1)})L\bigl(S(a_{(2)})\otimes a_{(3)}\bigr)
		-L\bigl(a_{(1)}S(a_{(2)})\otimes a_{(3)}\bigr)\\
	&\quad+L\bigl(a_{(1)}\otimes S(a_{(2)})a_{(3)}\bigr)
		-L\bigl(a_{(1)}\otimes S(a_{(2)})\bigr)\delta(a_{(3)})\\
	&=L(S(a_{(1)})\otimes a_{(2)})-L(a_{(1)}\otimes S(a_{(2)}).
\end{split}
\end{equation}
\end{proof}

\begin{Def}
Let $\mathcal{B}$ be a Hopf algebra and $L$ generator of an additive deformation. Then we set
\begin{equation}
\sigma:=L\circ(\id\otimes S)\circ\Delta=L\circ(S\otimes\id)\circ\Delta.
\end{equation}
\end{Def}

We freely choose between the two possibilities for $\sigma$, but Lemma\,\ref{Le:sigma} will only be essential in the proofs of Theorem\,\ref{Sa:constant_antipodes} and Lemma\,\ref{Le:sigma_circ_S} in section\,\ref{sec:const-antip-cocomm}.

\begin{Le}\label{Le:sigma-commuting}
$\sigma$ is commuting, i.e.
\begin{equation}
(\sigma\otimes\id)\circ\Delta=(\id\otimes\sigma)\circ\Delta.
\end{equation}
\end{Le}

\begin{proof}
% First we observe that
% \begin{equation}\begin{split}
% L(a_{(1)}\otimes S(a_{(2)}))&=L(a_{(1)}\otimes S(a_{(4)}))a_{(2)}S(a_{(3)})\\
% &=L\star\mu (a_{(1)}\otimes S(a_{(2)}))\\
% &=\mu\star L(a_{(1)}\otimes S(a_{(2)}))\\
% &=L(a_{(2)}\otimes S(a_{(3)}))a_{(1)}S(a_{(4)}).
% \end{split}\end{equation}
%Now we calculate
We calculate
\begin{equation}\begin{split}
(\sigma\otimes\id)\circ\Delta(a)&=\sigma(a_{(1)})a_{(2)}\\
&=(L\circ (\id\otimes S)\circ\Delta )(a_{(1)})\,a_{(2)}\\
&=L(a_{(1)}\otimes S(a_{(2)}))a_{(3)}\\
%&=L(a_{(2)}\otimes S(a_{(3)}))a_{(1)}S(a_{(4)})a_{(5)}\\
&=L(a_{(1)}\otimes S(a_{(4)}))a_{(2)}S(a_{(3)})a_{(5)}\\
&=a_{(1)}S(a_{(4)})L(a_{(2)}\otimes S(a_{(3)}))a_{(5)}\\
&=a_{(1)}L(a_{(2)}\otimes S(a_{(3)}))\\
&=(\id\otimes\sigma )\circ\Delta(a),
\end{split}\end{equation}
where the fifth equality holds since $L$ is commuting.
\end{proof}

\begin{Le}\label{Le:L-sigma} The following equations hold:
\begin{itemize}
\item
$
L^{\star n}\circ(\id\otimes S)\circ\Delta=\sigma^{\star n}
$
\item
$
e_\star^{tL}\circ (\id\otimes S)\circ\Delta=e_\star^{t\sigma}
$
\end{itemize}
\end{Le}

\begin{proof}
We prove this by induction over $n$. For $n=0,1$ the proposition is clear. We calculate
\begin{equation}\begin{split}
L^{\star n+1}(a_{(1)}\otimes S(a_{(2)}))
&=L\star L^{\star n}(a_{(1)}\otimes S(a_{(2)}))\\
&=L(a_{(1)}\otimes S(a_{(4)}))L^{\star n}(a_{(2)}\otimes S(a_{(3)}))\\
&=L(a_{(1)}\otimes S(a_{(3)}))\sigma^{\star n}(a_{(2)})\\
% &=L(a_{(1)}\otimes S(a_{(2n+2)}))L(a_{(2)}\otimes S(a_{(3)}))\dots L(a_{(2n)}\otimes S(a_{(2n+1)}))\\
&=L(a_{(1)}\sigma^{\star n}(a_{(2)})\otimes S(a_{(3)}))\\
&=L(\sigma^{\star n}(a_{(1)})a_{(2)}\otimes S(a_{(3)}))\\
&=\sigma^{\star n}(a_{(1)})\sigma(a_{(2)})\\
&=\sigma^{\star n+1}(a),
\end{split}\end{equation}
where the fifth equality follows from
Lemma\,\ref{Le:sigma-commuting}. The second equation follows easily
now:
\begin{equation}
e_\star^{tL}\circ(\id\otimes S)\circ\Delta
=\sum_{n=0}^\infty\frac{t^n}{n!}L^{\star n}\circ(\id\otimes S)\circ\Delta
=\sum_{n=0}^\infty\frac{t^n}{n!}\sigma^{\star n}
=e_\star^{t\sigma}.
\end{equation}
\end{proof}

\begin{Sa}\label{Sa:general_formula_for_S_t}
Let $\mathcal{B}$ be a Hopf algebra and $L$ generator of an additive deformation. Then it is a Hopf deformation and the deformed antipodes are given by 
\begin{equation}
S_t=S\star e_\star^{-t\sigma}.
\end{equation}
\end{Sa}

\begin{proof}
We have to check \eqref{eq:deformed_antipode}. By Lemma\,\ref{Le:L-sigma}
\begin{equation}\begin{split}
\mu_t\circ (\id\otimes S_t)\circ\Delta(a)
&=e_\star^{tL}\star \mu (a_{(1)}\otimes S(a_{(2)})) e_\star^{-t\sigma}(a_{(3)})\\
&=e_\star^{tL}(a_{(1)}\otimes S(a_{(4)}))a_{(2)}S(a_{(3)})
e_\star^{-t\sigma}(a_{(5)})\\
&=e_\star^{tL}(a_{(1)}\otimes S(a_{(2)}))e_\star^{-t\sigma}(a_{(3)})\Eins\\
&=e_\star^{t\sigma}(a_{(1)})e_\star^{-t\sigma}(a_{(2)})\Eins\\
&=\delta(a)\Eins.
\end{split}\end{equation}
The second equality in \eqref{eq:deformed_antipode} follows in the same manner.
\end{proof}

\goodbreak
\subsection{Constant antipodes in the cocommutative case}
\label{sec:const-antip-cocomm}

\begin{Le}
We have
\begin{equation}
\partial\sigma=L+L\circ (S\otimes S)\circ\tau.
\end{equation}
\end{Le}

\begin{proof}
\begin{equation}
\begin{split}
\partial\sigma(a\otimes b)&=\delta(a)\sigma(b)-\sigma(ab)+\sigma(a)\delta(b)\\
&=\delta(a) L(S(b_{(1)})\otimes b_{(2)})
	-L(S(a_{(1)}b_{(1)})\otimes a_{(2)}b_{(2)})\\
&\quad+L(S(a_{(1)})\otimes a_{(2)})\delta(b)\\
&=\delta(a) L(S(b_{(1)})\otimes b_{(2)})
	-L(S(b_{(1)})S(a_{(1)})\otimes a_{(2)}b_{(2)})\\
&\quad+L(S(a_{(1)})\otimes a_{(2)})\delta(b)\\
&=L(S(b)\otimes S(a))-L(S(a_{(1)})\otimes a_{(2)}b)+\delta(b) L(S(a_{(1)})\otimes a_{(2)})\\
&=L(S(b)\otimes S(a))+L(a\otimes b),
\intertext{
where in the fourth equality we used}
0&=\partial L(S(b_{(1)})\otimes S(a_{(1)})\otimes a_{(2)}b_{(2)})\\
&=\delta (b_{(1)})L(S(a_{(1)})\otimes a_{(2)}b_{(2)})
-L(S(b_{(1)})S(a_{(1)})\otimes a_{(2)}b_{(2)})\\
&\quad +L(S(b_{(1)})\otimes S(a_{(1)})a_{(2)}b_{(2)})
-L(S(b_{(1)})\otimes S(a_{(1)}))\delta (a_{(2)}b_{(2)})\\
&=L(S(a_{(1)})\otimes a_{(2)}b)
-L(S(b_{(1)})S(a_{(1)})\otimes a_{(2)}b_{(2)})\\
&\quad +\delta(a)L (S(b_{(1)})\otimes b_{(2)})
-L(S(b)\otimes S(a))
\intertext{
and in the last equality}
0&=\partial L(S(a_{(1)})\otimes a_{(2)}\otimes b)\\
&=\delta (a_{(1)})L(a_{(2)}\otimes b)
-L(S(a_{(1)})a_{(2)}\otimes b)\\
&\quad +L(S(a_{(1)})\otimes a_{(2)}b)
-L(S(a_{(1)})\otimes a_{(2)})\delta (b)\\
&=L(a\otimes b)
-\delta(a) L(\Eins\otimes b)\\
&\quad +L (S(a_{(1)})\otimes a_{(2)}b)
-L(S(a_{(1)})\otimes a_{(2)})\delta (b).
\end{split}\end{equation}
\end{proof}

\begin{Sa}\label{Sa:constant_antipodes}
Let $\mathcal{B}$ be a Hopf algebra, $L$ generator of an additive deformation. \\
If $\sigma=\sigma\circ S$
\begin{equation}
\widetilde{L}=L-\frac12\partial\sigma
\end{equation}
is the generator of a Hopf deformation with constant antipodes, i.e.
\begin{equation}
\widetilde{\mu}_t\circ(S\otimes\id)\circ\Delta
=\Eins\delta
=\widetilde{\mu}_t\circ (\id\otimes S)\circ\Delta.
\end{equation}
\end{Sa}

\begin{proof}
We can write
\begin{equation}
L=\underbrace{\frac{1}{2}(L+L\circ(S\otimes S)\circ\tau)}_{:=L_1}
+\underbrace{\frac{1}{2}(L-L\circ(S\otimes S)\circ\tau)}_{:=L_2}
\end{equation}
Then we have $L_1=\partial\frac{\sigma}{2}$ and 
$\sigma_2=L_2\circ(S\otimes\id)\circ\Delta=0$, since
\begin{equation}\begin{split}
L\circ (S\otimes S)\circ\tau\circ(S\otimes\id)\circ\Delta
&=L\circ(\id\otimes S)\circ (S\otimes S)\circ\tau\circ\Delta\\
&=L\circ(S\otimes\id)\circ\Delta\circ S\\
&=\sigma\circ S=\sigma,
\end{split}\end{equation}
where we made essential use of Lemma\,\ref{Le:sigma}.
\end{proof}

\begin{Le}\label{Le:sigma_circ_S}
If $\mathcal{B}$ is cocommutative, we have
\begin{equation}
\sigma=\sigma\circ S.
\end{equation}
\end{Le}

\begin{proof}
We calculate
\begin{equation}\begin{split}
\sigma\circ S&=L\circ(S\otimes\id)\circ\Delta\circ S\\
&=L\circ (S\otimes\id )\circ (S\otimes S)\circ\tau\circ\Delta\\
&=L\circ (S^2\otimes S)\circ\tau\circ\Delta\\
&=L\circ (\id\otimes S)\circ\Delta\\
&=\sigma,
\end{split}\end{equation}
due to Lemma\,\ref{Le:sigma}.
\end{proof}

So when deforming a cocommutative Hopf algebra one can always find an equivalent deformation, such that $S_t=S$ for all $t\in\mathbb{R}$.

\section{Examples}

\begin{Bsp}
In this example we realize the algebra of the quantum harmonic oscillator as the essentially only nontrivial additive deformation of the $*$-Hopf algebra of polynomials in adjoint commuting variables $\mathbb{C}\eklammer{x,x^*}$ with comultiplication and counit defined via
\begin{equation}
\Delta(x^\epsilon)=x^\epsilon\otimes\Eins+\Eins\otimes x^\epsilon
\qquad\text{and}\qquad
\delta(x^\epsilon)=0,
\end{equation} 
where $\epsilon\in\gklammer{1,*}$.

\begin{Beh}
Let $\mathcal{L}$ be an abelian Lie algebra, i.e.\ $\eklammer{a,b}=0\ \forall a,b\in\mathcal{L}$ and consider the universal enveloping Hopf algebra $U(\mathcal{L})$. In the case where $\mathcal{L}$ is of finite dimension $n$ this is just the polynomial algebra in $n$ commuting indeterminates.
For two additive deformations $\mu_t^{(1)},\mu_t^{(2)}$ of $U(\mathcal{L})$ with generators $L_1,L_2$ the following statements are equivalent:
\begin{enumerate}
\item $L_1-L_2$ is a coboundary i.e.\ the two deformations differ by a trivial deformation
\item $\mu_t^{(1)} (a\otimes b-b\otimes a)=\mu_t^{(2)} (a\otimes b-b\otimes a)$ for all $a,b\in\mathcal{L},t\in\mathbb{R}$
\item $L_1 (a\otimes b-b\otimes a)=L_2 (a\otimes b-b\otimes a)$ for all $a,b\in\mathcal{L},t\in\mathbb{R}$
\end{enumerate}
\end{Beh}

\begin{proof}
For any additive deformation of $U(\mathcal{L})$ we have
\begin{equation}\begin{split}
\mu_t (a\otimes b)
&=\mu\star e_\star^{tL}(a\otimes b)\\
&=\mu\otimes e_\star^{tL}
(a \otimes b \otimes \Eins \otimes \Eins
+a \otimes \Eins \otimes \Eins \otimes b
+\Eins \otimes b \otimes a \otimes \Eins
+\Eins \otimes \Eins \otimes a \otimes b)\\
&=ab+tL(a\otimes b)\Eins
\end{split}\end{equation}
as $L$ is normalized. From this the equivalence of 2 and 3 follows directly and to show that 1 is equivalent to 3 it suffices to show that $L$ is a coboundary iff $L(a\otimes b-b\otimes a)=0$ for all $a,b\in\mathcal{L}$, since we set $L=L_1-L_2$.

So let $L$ be a coboundary, i.e.\ $L=\partial\psi$. It follows that
\begin{equation}
L(a\otimes b-b\otimes a)
=-\psi(ab-ba)=0,
\end{equation}
since $\mathcal{L}$ is abelian.

Now let $L(a\otimes b-b\otimes a)=0$ for all $a,b\in\mathcal{L}$. Choose a basis of $\mathcal{L}$ and introduce any ordering on this bases. Then define
\begin{equation}
\psi (a_1\dots a_n):=
\begin{cases}
L(a_1\dots a_{n-1}\otimes a_n)&\text{ with $a_1\leq\dots\leq a_n$ if $n\geq 2$}\\
0&\text{ else.}
\end{cases}
\end{equation}
We write $\widetilde{L}=L+\partial\psi$ and 
$\widetilde{\mu}_t=\mu\star e_\star^{t\widetilde{L}}$. 
Now an easy induction on $n$ shows that $\widetilde{\mu}_t^{(n)}(a_1\otimes\dots \otimes a_n)=a_1\dots a_n$ for $a_1\leq\dots\leq a_n$. But from the equivalence of 2 and 3 we know that $\widetilde{\mu}_t$ is commutative so we get $\widetilde{\mu}_t=\mu$ for all $t\in\mathbb{R}$. So $\widetilde{L}=L+\partial\psi=0$ and $L$ is a coboundary.
\end{proof}

% \begin{Beh}
% Let $\mathcal{L}$ be an abelian Liealgebra, i.e.\ $\eklammer{a,b}=0\ \forall a,b\in\mathcal{L}$. Then two additive deformations of the universal enveloping Hopfalgebra $U(\mathcal{L})$ with generators $L_1,L_2$ are equivalent iff 
% \begin{equation}
% L_1(a\otimes b-b\otimes a)=L_2(a_\otimes b-b\otimes a)\qquad\forall a,b\in\mathcal{L}
% \end{equation}
% \end{Beh}

% \begin{proof}
% It clearly suffices to show that $L$ is generator of a trivial deformation iff $L(a\otimes b-b\otimes a)=0$, as $L_1,L_2$ generate equivalent deformations iff $L_1-L_2$ is trivial.

% Let $L$ be the generator of a trivial deformation. Then $L=\partial\psi$ for a normalized linear functional $\psi$. So we have
% \begin{equation}
% L(a\otimes b-b\otimes a)
% =-\psi(ab-ba)=0,
% \end{equation}
% since $\mathcal{L}$ is abelian.

% Now let $L(a\otimes b-b\otimes a)=0$ for all $a,b\in\mathcal{L}$. Introduce any ordering of the elements of $\mathcal{L}$ and define
% \begin{equation}
% \psi (a_1\dots a_n):=
% \begin{cases}
% L(a_1\dots a_{n-1}\otimes a_n)&\text{ with $a_1\leq\dots\leq a_n$ if $n\geq 2$}\\
% 0&\text{ else.}
% \end{cases}
% \end{equation}
% We write $\widetilde{L}=L-\partial\psi$ and 
% $\widetilde{\mu}_t=\mu\star e_\star^{t\widetilde{L}}$. 
% Now $\widetilde{\mu}_t$ coincides with $\mu$ on ordered expressions.  
% \end{proof}

It follows that a deformation of $\mathbb{C}\eklammer{x,x^*}$ is determined up to a trivial deformation by the value of $L(x\otimes x^*-x^*\otimes y)=\mu_1 (x\otimes x^*-x^*\otimes x)$. In case of a $*$-deformation $L$ must be hermitian, so this is a real number. Choosing different constants here corresponds to a rescaling of the deformation parameter $t$ so we assume $L(x\otimes x^*-x^*\otimes x)=1$. There is also a canonical representative for the cohomology class of the generator for which the antipodes are constant. Choosing $L(x\otimes x^*)=-L(x^*\otimes x)=1/2$ one gets $\sigma=0$.

One gets a well defined $*$-algebra isomorphism from the algebra generated by $a,a^\dagger$ and $\Eins$ with the relation $aa^\dagger-a^\dagger a=\Eins$ to the deformation of the polynomial algebra $(\mathbb{C}\eklammer{x,x^*},\mu_1)$ by setting $\Phi(a)=x$ and $\Phi(a^\dagger)=x^*$. In this sense the quantum harmonic oscillator algebra is the only nontrivial additive deformation of the polynomial algebra in two commuting adjoint variables.
\end{Bsp}

In the last three examples we take as Hopf algebra the group algebra $\mathbb{C}G$ over a group $G$. We identify linear functionals on $\mathbb{C}G^k$ with functions on $G^k$ for $k\in\mathbb{N}$.
For grouplike $a,b\in\mathcal{B}$ we have
\begin{equation}
\mu_t (a \otimes b) = e^{tL(a\otimes b)}ab.
\end{equation}

\begin{Bsp}\label{ex:Z}
We saw that in the cocommutative case it is possible to split an additive deformation into a trivial part and a part that corresponds to constant antipodes. But it is still possible that the part with constant antipodes is trivial as this example shows. Consider the $2$-coboundary defined by
\begin{equation}
L(m,n)=m^2n+mn^2
\end{equation} 
on the group algebra of $\mathbb{Z}$. In the following group elements of $\mathbb{Z}$ are denoted $(k)$ to avoid confusion with the complex number $k$.
This is a coboundary, since $L=\partial\psi$ where
\begin{equation}
\psi(k)=-\frac{1}{3}k^3
\end{equation}
We also see that $L(0,0)=0$ and $L$ is commuting, so $L\in B^{(\mathbf{NC})}$. Therefore it generates a trivial deformation. The deformation is nonconstant, since
\begin{equation}
\mu_t((1)\otimes(1))=e^{L((1),(1))}(2)=2(2)\ \neq\ (2)=\mu((1)\otimes(1)).
\end{equation}
We calculate
\begin{equation}
\sigma(k)=L((k),(-k))=-k^3+k^3=0
\end{equation}
for all $k\in\mathbb{Z}$, so the antipodes are constant. Since the deformation is trivial we can calculate the $\Phi_t$.
\begin{equation}
\Phi_t(k)=e^{t\psi(k)}k=e^{tk^3}k
\end{equation}
The second way for calculating the $S_t$ yields
\begin{equation}
S_t(k)=\Phi_t\circ S\circ\Phi_t(k)=e^{tk^3}\Phi_t(-k)=e^{tk^3}e^{-tk^3}(0)=(0).
\end{equation}
So in this situation we have $S\circ\Phi_t=\Phi_{-t}\circ S$.
\end{Bsp}

\begin{Bsp}
% We would like to see how Theorem\,\ref{Sa:constant_antipodes} can be find a simple represantitive for an equivalence class of additive deformations.
On $\mathbb{Z}^d$ every $d\times d$-matrix $A$ with complex entries defines a 2-cocycle $L$ via
\begin{equation}
L(\underline{k},\underline{l}):=\underline{k}A\underline{l}^t
\end{equation}
for $\underline{k},\underline{l}\in\mathbb{Z}^d$, since the functions $((k_1,\dots,k_d),(l_1,\dots,l_d))\mapsto k_il_j$ define cocycles for $i,j=1,\dots,d$, as is easily checked. These cocycles are of course normalized and commuting, so they are generators of additive deformations on a cocommutative Hopf algebra. $L$ is hermitian iff $A$ is hermitian. We want to apply 
Theorem\,\ref{Sa:constant_antipodes}, so we calculate
\begin{equation}
\sigma (\underline{k})=L(\underline{k},-\underline{k})=-\underline{k}A\underline{k}^t
\end{equation}
and
\begin{equation}\begin{split}
\partial\frac{\sigma}{2} (\underline{k},\underline{l})
&=\frac12 (-\underline{k}A\underline{k}^t+(\underline{k}+\underline{l})A(\underline{k}+\underline{l})^t
-\underline{l}A\underline{l}^t)\\
&=\frac12 (\underline{k}A\underline{l}^t+\underline{l}A\underline{k}^t)\\
&=\underline{k}\frac{A+A^t}{2}\underline{l}^t
\end{split}\end{equation}
which gives
\begin{equation}
\widetilde{L}(\underline{k},\underline{l})
=(L-\frac12\partial\sigma)(\underline{k},\underline{l}) 
=\underline{k}\frac{A-A^t}{2}\underline{l}^t.
\end{equation}
So every such cocycle is equivalent to one which comes from an antisymmetric matrix.
\end{Bsp}

\begin{Bsp}
Let $G$ be a group. then $\mathbb{C}G$ can be turned into a $*$-Hopf algebra in a natural way by extending the map $*:g\mapsto g^{-1}$ antilinearly to the whole of $\mathbb{C}G$. On the group elements the involution $*$ coincides with the antipode $S$. Now let $L$ be a generator of an additive $*$-deformation, i.e.\ $L$ is a normalized hermitian 2-cocycle. Then
\begin{equation}\begin{split}
\partial\frac{\sigma}{2}(g,h)
&=(L+L\circ (S\otimes S)\circ\tau )(g,h)\\
&=\frac{1}{2}(L(g,h)+L(h^*,g^*))\\
&=\frac{1}{2}(L(g,h)+\overline{L(g,h)})\\
&=\re L(g,h)
\end{split}\end{equation}
and consequently
\begin{equation}\begin{split}
\widetilde{L}(g,h)
&=L-\frac12\partial\sigma (g,h)\\
&=\im L(g,h).
\end{split}\end{equation}
So one has to consider only the case where $L$ is purely imaginary on the group elements.
\end{Bsp}

\begin{acknowledgements}
I would like to thank Prof.\ Michael Schürmann for many helpful discussions concerning this work. I also thank Karin Bokelmann and Stefan Voss for proofreading.
\end{acknowledgements}

%%% Local Variables: 
%%% mode: latex
%%% TeX-master: "Additive_Deformations"
%%% End: 

\bibliographystyle{myalphaurl-sortbyauthor}
\linespread{1.25}

\addcontentsline{toc}{chapter}{Literaturverzeichnis}
\bibliography{maltebib}

\end{document}